# COUNTEREXAMPLES TO A
# CONJECTURE OF LEMMERMEYER

NIGEL BOSTON AND CHARLES LEEDHAM-GREEN

January 7, 1998

**Abstract.** We produce infinitely many finite 2-groups that do not embed with index 2 in any group generated by involutions. This disproves a conjecture of Lemmermeyer and restricts the possible Galois groups of unramified 2-extensions, Galois over **Q**, of quadratic number fields.

## 1. Introduction.

Let $K$ be a quadratic number field, $H$ be its Hilbert class field, and $H_2$ its Hilbert 2-class field (in other words, the maximal everywhere unramified abelian and 2-abelian extensions of $K$ respectively). There has been much work recently regarding the possible structure of $\text{Gal}(H/K)$ and $\text{Gal}(H_2/K)$ (which are isomorphic to the ideal class group of $K$ and its 2-Sylow subgroup respectively). For instance, the conjectures of Cohen and Lenstra [2] suggest that every abelian group (respectively abelian 2-group) should appear as $\text{Gal}(H/K)$ (respectively $\text{Gal}(H_2/K)$) for a suitable choice of $K$.

A natural next question is to ask for nonabelian generalisations of this. For instance, one might ask what finite groups (respectively finite 2-groups) occur as the Galois group over $K$ of extensions obtained by iterating the construction of $H$ (respectively $H_2$). Note that such extensions are Galois over **Q**, not just over $K$. If we consider *any* unramified extension $L$ of $K$ with Galois group $G$ such that $L$ is Galois over **Q**, then $G$ has the following property:

($*$) $G$ embeds as a subgroup of index 2 in a group generated by elements of order 2 ("involutions").

This follows since $\text{Gal}(L/\mathbf{Q})$ is generated by its inertia subgroups, all of which have order at most 2, since $L/K$ is everywhere unramified. So a first step towards the question of the last paragraph is to ask whether every finite group (or in particular every finite 2-group) satisfies ($*$). Franz Lemmermeyer conjectured that every finite 2-group satisfies ($*$) (see, for example, http://www.rzuser.uni-heidelberg.de/~hb3/unsol.ps). In this paper we show that in fact there exist finite 2-groups that do not satisfy ($*$). They cannot, therefore, occur as Galois groups of iterated 2-class fields of quadratic fields and any generalisations of the Cohen-Lenstra heuristics will have to take this into account.

1991 *Mathematics Subject Classification*. 11R37, 20D15.
The first author was partially supported by NSF grant DMS 96-22590.

**2. The Counterexamples.**

It is straightforward to find the smallest counterexamples by a systematic check using the MAGMA databases [1] of groups of order dividing 256. It turns out that all groups of order dividing 32 satisfy $(*)$. Of the 267 groups of order 64, exactly 2 fail to satisfy $(*)$. These are numbers 28 and 46 in the database. They can be given by respective compact presentations:

$G_1 = \langle a, b | a^{16}, b^4, [b,a] = a^4 \rangle$,
$G_2 = \langle a, b | a^{16}, b^4, [b,a] = a^{-2} \rangle$.

The method employed is to calculate, for each of the 2328 groups of order 128, the subgroup generated by all involutions. Those groups for which this subgroup is the whole group are saved and their maximal subgroups computed. This produced a list of 7007 groups of order 64. These fall into 265 isomorphism classes. Checking these against the 267 groups of order 64 yields the 2 isomorphism classes missed, namely $G_1$ and $G_2$ above.

Similarly one can work with the maximal subgroups of groups of order 256 with $\leq$ 4 generators, generated by involutions. Of the 162 groups of order 128 with precisely 2 generators, exactly 11 fail to satisfy $(*)$, namely numbers $45, 87, 88, 98, 100, 101, 130, 144, 145, 146$, and 158 in the database. Of the 833 groups of order 128 with precisely 3 generators, exactly 7 fail to satisfy $(*)$, namely numbers $836, 841, 842, 886, 887, 912$, and 971 in the database. It is not computationally feasible to search for all the groups of order 128 with precisely 4 generators failing $(*)$.

**Lemma 1** Suppose that $G$ is a finite 2-group and $N$ a characteristic subgroup. If $G/N$ fails $(*)$, then so does $G$.

*Proof* Suppose, on the contrary, that $G$ embeds with index 2 in a group $\Gamma$ generated by involutions. Then $N$ is normal in $\Gamma$ and $G/N$ embeds with index 2 in $\Gamma/N$, which is generated by involutions.

Using this idea, we obtain at least 2 further examples of groups failing $(*)$, these of order 256. Numbers 45 and 100 in the database of groups of order 128 have nuclear rank 1 (in the terminology of [3]) and each have 1 (terminal) descendant. All other groups mentioned above have nuclear rank 0.

It is interesting to note that 4 of the 18 known counterexamples of order 128 (namely numbers $45, 100, 841, 886$) have quotients isomorphic to our counterexamples of order 64, but that in each case the kernel is not characteristic.

It is possible to give human proofs that the groups above fail to satisfy $(*)$. For instance, the following argument shows that there exist infinitely many counterexamples to Lemmermeyer's conjecture.

**Lemma 2** Let $G = G(m,n) := \langle a, b | a^m, b^n, [b,a] = a^4 \rangle$, where $m$ and $n$ are powers of 2, $m \geq 16, n \geq 4, m \leq 4n$. Then $G$ has order $mn$, $G' = \langle a^4 \rangle$ is cyclic of order $m/4$, and $G/G'$ is a product of cyclic groups of order 4 and $n$.

*Proof* Immediate.

*Remark* $G(16,4)$ is $G_1$ above. $G(16,8)$ is group number 45 in the database of groups of order 128. $G(32,8)$ is the descendant of order 256 of that group, mentioned above.

**Lemma 3** Let $G = G(m,n)$, as above. Let $N = \langle a^4 \rangle$. Then $N$ is a characteristic subgroup of $G$, and $G/N$ is abelian, but no automorphism of $G$ induces the inverting automorphism $g \mapsto g^{-1}$ on $G/N$.

*Proof* Suppose that such an automorphism $\alpha$ exists. Then $\alpha(a) = a^{-1}a^{4r}, \alpha(b) = b^{-1}a^{4s}$, for some integers $r, s$. We calculate $\alpha(a^4)$ in two ways.

First, $\alpha(a^4) = \alpha(a)^4 = (a^{4r-1})^4$.

Second, note that since in any group $[xy, z] = [x,z]^y[y,z]$, it follows that for any integers $u, v, w$, we have $[b^u a^v, a^w] = [b,a]^{uw} = a^{4uw}$. So $\alpha(a^4) = \alpha([b,a]) = [b^{-1}a^{4s}, a^{4r-1}] = a^{-4(4r-1)}$.

Equating the two expressions gives that the order of $a$ divides $32r - 8$, a contradiction.

**Lemma 4** Let $A$ be an abelian 2-group, and let $\alpha$ be an automorphism of $A$ of order 2. Then the set of elements of $A$ inverted by $\alpha$ forms a subgroup of $A$.

*Proof* Trivial.

**Lemma 5** Let $A$ be an abelian 2-group, and let $B$ be a group containing $A$ as a subgroup of index 2. If $B$ is generated by involutions, then $A$ is the product of the subgroup $\Omega(A)$ (consisting of the elements of $A$ of order dividing 2) with the subgroup consisting of those elements of $A$ that are inverted by an element of $B \setminus A$.

*Proof* Let $a \in A$. Write $a = c_1 c_2 \ldots c_t$, where the $c_i$ are involutions. Replacing $c_i$ by a conjugate as appropriate, we may assume that the $c_i$ in the above expression that lie in $A$ are $c_1, \ldots, c_r$ for some $r \geq 0$. Since $a \in A$ it follows that $t - r$ is even, but the product of two involutions is inverted by either of the involutions, and so the result follows.

**Theorem** Let $G = G(m,n)$, as defined in Lemma 2. Then $G$ fails to satisfy $(*)$. Thus, there are infinitely many finite 2-groups that fail $(*)$.

*Proof* Suppose that $G$ embeds with index 2 in a group $\Gamma$ generated by involutions. Then $\Gamma$ contains $G$ as a normal subgroup, and hence $N$ (as in Lemma 3) is a normal subgroup of $\Gamma$. By Lemma 5, since $\Gamma/N$ is generated by involutions, $G/N$ is the product of $\Omega(G/N)$ with the subgroup of $G/N$ inverted by an element of $\Gamma \setminus G$. But $\Omega(G/N)$ is in the Frattini subgroup of $G/N$, so $\Gamma$ contains an involution $\alpha$ that inverts every element of $G/N$; but this is impossible by Lemma 3.

**References**

1. W.Bosma and J.J.Cannon, *Handbook of Magma Functions* (1996), Sydney.
2. H.Cohen and H.W.Lenstra, Jr., *Heuristics on class groups of number fields*, Number theory, Noordwijkerhout 1983, Lecture Notes in Math. 1068, Springer, Berlin-New York, 1984, pp. 33-62.
3. E.A.O'Brien, *The p-group generation algorithm*, J. Symbolic Computation **9** (1990), no. 5-6, 677-698.

Department of Mathematics, University of Illinois, Urbana, Illinois 61801
*E-mail address*: `boston@math.uiuc.edu`

School of Mathematical Sciences, Queen Mary and Westfield College, University of London, Mile End Road, London E1 4NS
*E-mail address*: `C.R.Leedham-Green@qmw.ac.uk`